\documentclass[a4paper,10pt]{article}

\usepackage[utf8]{inputenc}
\usepackage{amsmath,amssymb}
\newtheorem{thm}{Theorem}[section]

\newtheorem{pro}[thm]{Proposition}
\newtheorem{cor}[thm]{Corollary}
\usepackage[all]{xy}

\title{$T$-stability for Higgs sheaves over compact complex manifolds}
\author{S. A. H. Cardona\footnote{Electronic addresses: s.a.h.cardona@gmail.com, andres.holguin@cimat.mx}\\CIMAT A.C. - Via Jalisco S/N - 36240, Gto. - M\'exico}

\begin{document}

\maketitle

\begin{abstract}
We introduce the notion of $T$-stability for torsion-free Higgs sheaves as a natural generalization of the notion of $T$-stability 
for torsion-free coherent sheaves over compact complex manifolds. We prove similar properties to the 
classical ones for Higgs sheaves. In particular, we show that only saturated flags of torsion-free Higgs sheaves 
are important in the definition of $T$-stability. Using this, we show that this notion is preserved under 
dualization and tensor product with an arbitrary Higgs line bundle. Then, we prove that for a torsion-free 
Higgs sheaf over a compact K\"ahler manifold, $\omega$-stability implies $T$-stabilty. As a consequence of this we 
obtain the $T$-semistability of any reflexive Higgs sheaf with an admissible Hermitian-Yang-Mills metric. Finally, we prove that 
$T$-stability implies $\omega$-stability if, as in the classical case, some additional requirements on the base 
manifold are assumed. In that case, we obtain the existence of admissible Hermitian-Yang-Mills metrics on any 
$T$-stable reflexive sheaf. \\


\noindent{\it MS Classification}: 53C07, 53C55, 32C15.

\end{abstract}

\section{Introduction}

The notion of $T$-stability was introduced by Bogomolov \cite{Bogomolov} in the case of coherent sheaves over 
projective algebraic manifolds, and it was studied latter by Kobayashi in \cite{Kobayashi 0} and \cite{Kobayashi} for 
coherent sheaves over compact complex manifolds. In the K\"ahler case, the $T$-stability was related to 
Mumford-Takemoto stability (also called $\omega$-stability, where $\omega$ denotes the K\"ahler form of the base manifold). 
To be precise, it was shown by Bogomolov and Kobayashi that a $\omega$-stable (resp. $\omega$-semistable) 
torsion-free coherent sheaf over a compact K\"ahler manifold $X$, was $T$-stable (resp. $T$-semistable). They proved also the converse 
result if $H^{1,1}(X,{\mathbb C})$ was one dimensional, or if $\omega$ represented an integral class (so that $X$ was projective algebraic) and 
${\rm Pic}(X)/{\rm Pic}^{0}(X)={\mathbb Z}$, where 
here ${\rm Pic}^{0}(X)$ denotes the subgroup of the Picard group ${\rm Pic}(X)$ consisting of holomorphic 
line bundles with vanishing first Chern class. In proving the connection between these two concepts of stability, 
it was important to consider a classical vanishing theorem for holomorphic line bundles. As we will see, the same 
result is important also to connect these two notions of stability for Higgs sheaves.\\

Now, vanishing theorems are important in Complex Geometry. Indeed, some of these results, first 
proved by Bochner and Yano \cite{Yano-Bochner}, have been used by Kobayashi \cite{Kobayashi} to prove one 
direction of the classical Hitchin-Kobayashi correspondence for holomorphic vector bundles over compact 
K\"ahler manifolds. As it is well known, this correspondence establishes an equivalence between the notion of 
$\omega$-polystabilty and the existence of Hermitian-Einstein metrics for such bundles. Kobayashi also proved 
that a holomorphic vector bundle admitting an approximate Hermitian-Einstein metric was $\omega$-semistable. 
As a consequence of this, it was followed that a holomorphic vector bundle over a compact K\"ahler manifold admitting a 
Hermitian-Einstein metric (resp. an approximate Hermitian-Einstein metric) was necessarily $T$-stable 
(resp. $T$-semistable). The Hitchin-Kobayashi correspondence has been extended to reflexive sheaves by Bando 
and Siu \cite{Bando-Siu} by introducing the notion of admissible metric on a sheaf.\\

On the other hand, Higgs bundles and Higgs sheaves were introduced by Hitchin \cite{Hitchin} and Simpson 
\cite{Simpson}, \cite{Simpson 2} and they also introduced the corresponding notion 
of Mumford-Takemoto stability for these objects. As it is well known, several results on holomorphic vector bundles and coherent 
sheaves can be extended to Higgs bundles and Higgs sheaves. In particular, Vanishing theorems for Higgs 
bundles have been recently studied in \cite{Cardona 3}, and Simpson proved in \cite{Simpson} a Hitchin-Kobayashi 
correspondence for Higgs bundles over compact K\"ahler manifolds, i.e., an equivalence between the notion of 
Mumford-Takemoto polystability and the existence of Hermitian-Yang-Mills metrics 
(henceforth usually abbreviated $HYM$-metric). Now, Bruzzo and Gra\~na Otero \cite{Bruzzo-Granha} proved that if 
a Higgs bundle admits an approximate Hermitian-Yang-Mills (henceforth abbreviated $apHYM$-metric), it is
necessarily semistable in the sense of Mumford-Takemoto. Following the ideas of Bando and Siu \cite{Bando-Siu}, 
Biswas and Schumacher proved in \cite{Biswas-Schumacher} that a reflexive Higgs sheaf over a compact  K\"ahler 
manifold is $\omega$-polystable if and only if it has an admissible Hermitian-Yang-Mills metric, which is indeed a 
Hitchin-Kobayashi correspondence for Higgs sheaves. \\

This article is organized as follows, in the second section we review some basic definitions concerning Higgs 
sheaves and using a classical isomorphism involving determinants, we show that the determinant bundle of 
certain quotient Higgs sheaves is indeed a Higgs line bundle. In the final part of the second section, we review 
some results on reflexive Higgs sheaves and Higgs bundles and we rewrite a classical vanishing theorem of 
holomorphic line bundles over compact K\"ahler manifolds in the context of Higgs line bundles. In the third 
section, we introduce the notion of a weighted flag for a torsion-free Higgs sheaf, and using this we define 
the $T$-stability for torsion-free Higgs sheaves, as a natural extension of the notion of $T$-stability for 
torsion-free coherent sheaves introduced in \cite{Kobayashi}. Then, we prove some basic properties that are indeed 
extensions of classical results; in particular, we prove that in the definition of $T$-stability it is enough 
to consider saturated flags (i.e., flags in which the quotients between the Higgs sheaf and all Higgs subsheaves of 
the flag are torsion-free). We prove also that $T$-stability is preserved under tensor products with Higgs 
line bundles and dualizations. Finally, in a four section we prove that in the K\"ahler case, a Mumford-Takemoto stable (resp. semistable) torsion-free Higgs sheaf 
is also $T$-stable (resp. $T$-semistable) and hence, as a consequence of the main result of Biswas and Schumacher 
in \cite{Biswas-Schumacher}, we obtain that any reflexive Higgs sheaf with an admissible $HYM$-metric is 
necessarily $T$-semistable and it is in general a direct sum of $T$-stable Hermitian-Yang-Mills Higgs sheaves with 
equal slope; and that any locally free Higgs sheaf admitting a $HYM$-metric (resp. $apHYM$-metric) is $T$-stable 
(resp. $T$-semistable). At the end, we show that if either $H^{1,1}(X,{\mathbb C})$ is one dimensional, or if 
$\omega$ represents an integral class and ${\rm Pic}(X)/{\rm Pic}^{0}(X)={\mathbb Z}$, the classical proof of 
Kobayashi for the converse implication (i.e., $T$-stability as a sufficient condition of $\omega$-stability) can 
be easily adapted to the Higgs case.

\section{Preliminaries}

We start with some basic definitions. Let $X$ be a compact complex manifold of complex dimension $n$ and let 
$\Omega_{X}^{1}$ be the cotangent sheaf to $X$. A Higgs sheaf over $X$ is a pair ${\mathfrak E}=(E,\phi)$ where 
$E$ is a coherent sheaf over $X$ and $\phi : E\rightarrow E\otimes\Omega_{X}^{1}$ is a morphism of ${\cal O}_{X}$-modules 
such that $\phi\wedge\phi : E\rightarrow E\otimes\Omega_{X}^{2}$ vanishes. The morphism $\phi$ is usually called 
the Higgs field. A section $s$ of $E$ is said to be a $\phi$-invariant section of ${\mathfrak E}$, if there exists 
a holomorphic 1-form $\lambda$ of $X$ such that $\phi(s)=s\otimes\lambda$. A Higgs sheaf ${\mathfrak E}$ is said 
to be torsion-free (resp. locally free, reflexive, normal, torsion) if the coherent sheaf $E$ is torsion-free 
(resp. locally free, reflexive, normal, torsion). The support of a Higgs sheaf is the support of the corresponding 
coherent sheaf, and hence ${\rm supp}(\mathfrak E)=\{x\in X; E_{x}\neq 0\}$. If ${\mathfrak T}=(T,\psi)$ is a 
torsion Higgs sheaf, from a classical result of Kobayashi \cite{Kobayashi}, we know that ${\rm det}\,T$ admits a 
nonzero holomorphic section and, in particular, if ${\rm supp}(\mathfrak T)$ has codimension at least two, 
${\rm det}\,{\mathfrak T}$ is a trivial holomorphic line bundle. \\

As it is well known \cite{Biswas-Schumacher}, on Higgs sheaves we can apply the same operations that we normaly 
apply to sheaves. For instance, the dual of a Higgs sheaf and its pullback are again Higgs sheaves, and tensor 
products and the direct sums of Higgs sheaves are Higgs sheaves. If ${\mathfrak E}$ is a Higgs sheaf we denote 
its dual by ${\mathfrak E}^{*}$, and if $f:Y\longrightarrow X$ is a map between compact complex manifolds, 
we denote its pullback by $f^{*}{\mathfrak E}$. If now ${\mathfrak E}_{1}$ and ${\mathfrak E}_{2}$ are Higgs 
sheaves, we denote its tensor product and direct sum by ${\mathfrak E}_{1}\otimes{\mathfrak E}_{2}$ and 
${\mathfrak E}_{1}\oplus{\mathfrak E}_{2}$ respectively. Now, a Higgs subsheaf 
${\mathfrak F}$ of ${\mathfrak E}$ is a subsheaf $F$ of $E$ such that $\phi(F)\subset F\otimes\Omega_{X}^{1}$, 
and hence the pair ${\mathfrak F}=(F,\phi|_{F})$ becomes itself a Higgs sheaf. A morphism $f:{\mathfrak E}_{1}\longrightarrow{\mathfrak E}_{2}$ 
between two Higgs sheaves over $X$, is a morphism $f:E_{1}\longrightarrow E_{2}$ of the corresponding coherent 
sheaves such that the diagram
 \begin{displaymath}
 \xymatrix{
E_{1} \ar[r]^{\phi_{1}} \ar[d]^{f}    &     E_{1}\otimes\Omega_{X}^{1} \ar[d]^{f\otimes 1}   \\
E_{2} \ar[r]^{\phi_{2}}               &     E_{2}\otimes\Omega_{X}^{1}   \\
}
\end{displaymath}
is commutative. If ${\mathfrak E}=(E,\phi)$ is a Higgs sheaf over $X$, the natural morphism $\sigma:E\rightarrow E^{**}$ is 
a first example of a Higgs morphism $\sigma:{\mathfrak E}\rightarrow {\mathfrak E}^{**}$. The kernel and the image of 
Higgs morphisms are Higgs sheaves and the torsion subsheaf of a Higgs sheaf is again a Higgs sheaf 
(see \cite{Cardona 2} for details), these two results will be particularly important in the study of 
$T$-stability. An exact sequence of Higgs sheaves is an exact sequence of the corresponding coherent sheaves in 
which each morphism is a morphism of Higgs sheaves. \\

Let ${\mathfrak E}$ be a torsion-free Higgs sheaf of rank $r$, from a classical result (see \cite{Kobayashi} for 
details of this and what follows) we know that ${\rm det\,}E\cong\left(\bigwedge^{r}E\right)^{**}$ and hence, by 
using this isomorphism, it is possible to induce a Higgs field on the determinant bundle. As a consequence of this, 
we see that the determinant bundle of a torsion-free Higgs sheaf is a Higgs line bundle, we denote this bundle by 
${\rm det\,}{\mathfrak E}$. Clearly, from this definition we have canonically ${\rm det}\,{\mathfrak E}\cong\left(\bigwedge^{r}{\mathfrak E}\right)^{**}$ 
as an isomorphism of Higgs bundles. Now, let us consider the short exact sequence of Higgs sheaves 
 \begin{equation}
 \xymatrix{
0 \ar[r]  &  {\mathfrak F} \ar[r]  &   {\mathfrak E} \ar[r]  &   {\mathfrak G} \ar[r] &  0  
} \nonumber
\end{equation}
(also called a Higgs extension). If ${\mathfrak E}$ is torsion-free, then ${\mathfrak F}$ is torsion-free, but 
${\mathfrak G}$ may have torsion. In this case, we induce a Higgs morphism on ${\rm det}\,G$ using the Higgs 
fields of ${\rm det}\,{\mathfrak F}$ and ${\rm det}\,{\mathfrak E}$, and the isomorphism 
${\rm det}\,G\cong({\rm det}\,F)^{-1}\otimes{\rm det}\,E$. We denote by ${\rm det}\,{\mathfrak G}$ the Higgs line 
bundle defined by ${\rm det}\,G$ and this induced morphism. In this way we obtain 
${\rm det}\,{\mathfrak E}\cong{\rm det}\,{\mathfrak F}\otimes{\rm det}\,{\mathfrak G}$ as an isomorphism of Higgs 
bundles.\\

Suppose now that $X$ is a K\"ahler manifold with $\omega$ its K\"ahler form, then the first Chern class of 
$\mathfrak E$ is by definition the first Chern class of $E$, and hence following Kobayashi \cite{Kobayashi}, 
$c_{1}(\mathfrak E)=c_{1}({\rm det}\,E)$ and the degree of $\mathfrak E$ is given by 
\begin{equation}
 {\rm deg\,}{\mathfrak E} = \int_{X}c_{1}(\mathfrak E)\wedge\omega^{n-1}\,.  \label{deg}
\end{equation}
It is important to note that the degree defined by (\ref{deg}) depends on $\omega$ if the complex dimension of $X$ 
is greater than one. Now, if we denote the rank of $\mathfrak E$ by ${\rm rk\,}{\mathfrak E}$, and if this rank is 
positive, we introduce the quotient $\mu({\mathfrak E})={\rm deg\,}{\mathfrak E}/{\rm rk\,}{\mathfrak E}$, which 
is called the slope of the Higgs sheaf. A Higgs sheaf $\mathfrak E$ is said to be $\omega$-stable (resp. $\omega$-semistable), if it is torsion-free 
and for any Higgs subsheaf ${\mathfrak F}$ with $0<{\rm rk\,}{\mathfrak F}<{\rm rk\,}{\mathfrak E}$ 
we have the inequality $\mu({\mathfrak F})<\mu(\mathfrak E)$ (resp. $\le$). We say that a Higgs sheaf is $\omega$-polystable if it decomposes into a direct sum of two or more 
$\omega$-stable Higgs sheaves all these with the same slope.\\

This notion of stability was introduced by Hitchin \cite{Hitchin} and Simpson \cite{Simpson} as an analog of the 
Mumford-Takemoto stability for coherent sheaves \cite{Kobayashi}. However, it is important to note that the 
coherent sheaf $E$ associated to a torsion-free Higgs sheaf ${\mathfrak E}$, is $\omega$-stable (resp. $\omega$-semistable) 
if and only if for any proper nontrivial subsheaf $F$ of $E$, we have $\mu(F)<\mu(E)$ (resp. $\le$). Therefore, if 
${\mathfrak E}$ is $\omega$-stable (resp. $\omega$-semistable) in the classical sense, it is $\omega$-stable 
(resp. $\omega$-semistable) as a Higgs object, but the converse is not true in general (see \cite{Hitchin} for examples). 
In this sense, the notion of $\omega$-stability for Higgs sheaves is a generalization of the classical notion of 
$\omega$-stability for coherent sheaves.\\

As it is well known (see for instance \cite{Cardona 2} or \cite{Simpson}), for this notion of 
stability it is suffice to consider only Higgs subsheaves with torsion-free quotients (we will see in the next 
section that there exists a similar result for $T$-stability). As we said before, Biswas and Schumacher proved  
\cite{Biswas-Schumacher} the equivalence between $\omega$-stability and the existence of $HYM$-metrics for Higgs 
sheaves. To be precise, they proved the following result:
\begin{thm}\label{HK-corresp. sheaves}
Let ${\mathfrak E}$ be a reflexive Higgs sheaf over a compact K\"ahler manifold $X$ with K\"ahler form $\omega$. 
Then, there exists an admissible HYM-metric on ${\mathfrak E}$ if and only if it is $\omega$-polystable.
\end{thm}

A Higgs bundle is by definition a locally free Higgs sheaf. We say that a Higgs bundle is Hermitian flat if there 
exists a Hermitian metric $h$ on it, such that the Hitchin-Simpson connection ${\cal D}_{h}=D_{h}+\phi+\bar\phi_{h}$ is flat, 
i.e., if the Hitchin-Simpson curvature ${\cal R}_{h}={\cal D}_{h}\wedge{\cal D}_{h}$ vanishes. Now, following
\cite{Bruzzo-Granha} we know that
\begin{equation}
 {\cal R}_{h} = R_{h} + D'_{h}(\phi) + D''(\bar\phi_{h}) + [\phi,\bar\phi_{h}] \label{R Hitchin-Simpson}
\end{equation}
where $R_{h}$ is the Chern curvature, $D'_{h}$ and $D''$ are the holomorphic and anti-holomorphic parts of the 
Chern connection $D_{h}$ and the commutator is the usual abbreviation for $\phi\wedge\bar\phi_{h} + \bar\phi_{h}\wedge\phi$. 
Now, on the right hand side of (\ref{R Hitchin-Simpson}) the third term is the adjoint of the second term 
and since for Higgs line bundles the commutator is zero, a Higgs line bundle ${\mathfrak L}=(L,\phi)$ is Hermitian flat if 
and only if $R_{h}=0$ (i.e., $L$ is Hermitian flat in the classical sense) and the Higgs field satisfies 
$D'_{h}\phi=0$. Notice that in the case of Higgs line bundles, any Higgs morphism is in essence a holomorphic 
1-form, hence every holomorphic section of a Higgs line bundle is an invariant section. \\

On the other hand, in Complex Geometry there is a well known vanishing theorem for holomorphic line bundles 
depending on its degree \cite{Kobayashi}, since the degree of a Higgs bundle is the same degree of the 
corresponding vector bundle, this result can be applied to Higgs line bundles, and hence the classical 
vanishing theorem in the Higgs context becomes 
\begin{pro} \label{Prop. Line bundles}
 Let $\mathfrak L$ be a Higgs line bundle over a compact K\"ahler manifold $X$. Then \\
{\bf (i)} If ${\rm deg}\,{\mathfrak L}<0$, then ${\mathfrak L}$ admits no nonzero (invariant) holomorphic sections;\\
{\bf (ii)} If ${\rm deg}\,{\mathfrak L}=0$, then every nonzero (invariant) holomorphic section of ${\mathfrak L}$ has no zeros.
\end{pro}
Finally, as it is well known, in the case of Higgs bundles over compact K\"ahler manifolds, Simpson \cite{Simpson} 
proved a Hitchin-Kobayashi correspondence, and Bruzzo and Gra\~na Otero proved in \cite{Bruzzo-Granha} that Higgs 
bundles admitting $apHYM$-metrics are semistable in the sense of Mumford-Takemoto. The converse of this result has 
been proved in \cite{Cardona} in the one-dimensional case and by Li and Zhang \cite{Jiayu-Zhang} for compact 
K\"ahler manifolds of greater dimensions. These results can be summarized as follows:
\begin{thm}\label{HK-corresp.}
Let ${\mathfrak E}$ be a Higgs bundle over a compact K\"ahler manifold $X$ with K\"ahler form $\omega$. Then, 
there exists a HYM-metric (resp. apHYM-metric) on ${\mathfrak E}$ if and only if it is $\omega$-polystable 
(resp. $\omega$-semistable).
\end{thm}
Notice that, since for compact K\"ahler manifolds an admissible $HYM$-metric is just a $HYM$-metric, part 
of Theorem \ref{HK-corresp.} is indeed a particular case of Theorem \ref{HK-corresp. sheaves}. However, there is 
known a differential geometric analog of $\omega$-semistability only for bundles\footnote{Indeed, even 
in the classical case of reflexive sheaves, there is no yet an equivalence of $\omega$-semistability 
(see \cite{Bando-Siu} for more details).} and it is precisely the notion of $apHYM$-metric, a natural 
extension for Higgs bundles of an approximate Hermitian-Einstein structure for holomorphic vector bundles.

\section{$T$-stability}
In order to define the notion of $T$-stability, we need to define first the notion of a weighted flag in the 
Higgs case. Let ${\mathfrak E}$ be a torsion-free Higgs sheaf over a compact complex manifold $X$, a weighted 
flag of ${\mathfrak E}$ is a sequence of pairs ${\cal F}=\{({\mathfrak E}_{i},n_{i})\}_{i=1}^{k}$ consisting 
of Higgs subsheaves 
\begin{equation}
 {\mathfrak E}_{1}\subset{\mathfrak E}_{2}\subset\cdots\subset{\mathfrak E}_{k}\subset{\mathfrak E}  \nonumber
\end{equation}
together with positive integers $n_{1},n_{2},...,n_{k}$ and such that
\begin{equation}
 0<{\rm rk}\,{\mathfrak E}_{1}<{\rm rk}\,{\mathfrak E}_{2}<\cdots<
   {\rm rk}\,{\mathfrak E}_{k}<{\rm rk}\,{\mathfrak E}\,.     \nonumber
\end{equation}
Let $r_{i}={\rm rk}\,{\mathfrak E}_{i}$ and $r={\rm rk}\,{\mathfrak E}$. In analogy to the classical case, 
to each weighted flag ${\cal F}$ we associate the Higgs line bundle
\begin{equation}
 {\mathfrak T}_{\cal F}=\prod_{i=1}^{k}(({\rm det}\,{\mathfrak E}_{i})^{r}\otimes({\rm det}\,{\mathfrak E})^{-r_{i}})^{n_{i}}. \label{def T}
\end{equation}
We say that a weighted flag ${\cal F}$ is saturated if the quotients ${\mathfrak E}/{\mathfrak E}_{i}$, $i=1,2,...,k$, 
are all torsion-free. A torsion-free Higgs sheaf ${\mathfrak E}$ over $X$ is said to be $T$-stable 
(resp. $T$-semistable), if for every weighted flag ${\cal F}$ of ${\mathfrak E}$ and every Hermitian flat Higgs 
line bundle ${\mathfrak L}$ over $X$, the Higgs line bundle ${\mathfrak T}_{\cal F}\otimes{\mathfrak L}$ admits 
no nonzero holomorphic sections (resp. every nonzero holomorphic section of 
${\mathfrak T}_{\cal F}\otimes{\mathfrak L}$, if any, vanishes nowhere on $X$). Notice that since weighted flags 
for Higgs sheaves consist of Higgs subsheaves, if a Higgs sheaf is $T$-stable (resp. $T$-semistable) in the 
classical sense, i.e., as a coherent sheaf, it is also $T$-stable (resp. $T$-semistable) in the Higgs sense. 
However, as we will see in the next section the converse is not true in general.\\

From this definition of $T$-stability we have the following results, which are natural extensions to the Higgs 
case of classical results of Kobayashi \cite{Kobayashi}.
\begin{pro}\label{prop. saturated}
 Let ${\mathfrak E}$ be a torsion-free Higgs sheaf over a compact complex manifold $X$. Then it is $T$-stable 
 (resp. $T$-semistable) if and only if for every saturated flag ${\cal F}$ of ${\mathfrak E}$ and every Hermitian 
 flat Higgs line bundle $\mathfrak L$ over $X$, the bundle ${\mathfrak T}_{\cal F}\otimes{\mathfrak L}$ admits no 
 nonzero holomorphic sections (resp. every holomorphic section of ${\mathfrak T}_{\cal F}\otimes{\mathfrak L}$, if any, vanishes nowehere 
 on $X$).
\end{pro}
\noindent {\it Proof:} There is nothing to prove in one direction\footnote{If ${\mathfrak E}$ is $T$-stable 
(resp. $T$-semistable), the conditions on existence or not of holomorphic sections hold for any flag and any 
Hermitian flat Higgs line bundle; in particular this is true if the flag is saturated.}. Now, in 
order to prove the other direction, let us assume that such conditions on existence or not of holomorphic sections
are satisfied for any saturated flag and any Hermitian flat Higgs line bundle. \\

Let ${\cal F}'=\{({\mathfrak E}_{i},n_{i})\}_{i=1}^{k}$ be an arbitrary flag of ${\mathfrak E}$ and 
${\mathfrak L}$ a Hermitian flat Higgs line bundle. Let ${\mathfrak T}_{i}$ be the torsion of 
${\mathfrak E}/{\mathfrak E}_{i}$. Then, if we define $\tilde{\mathfrak E}_{i}$ as the kernel of the morphism 
${\mathfrak E}\rightarrow({\mathfrak E}/{\mathfrak E}_{i})/{\mathfrak T}_{i}$ we obtain the following 
commutative diagram
\begin{displaymath}
 \xymatrix{
         &                                &                            &            0   \ar[d]           &   \\
         &           0  \ar[d]            &                            &  {\mathfrak T}_{i} \ar[d]       &   \\
0 \ar[r] & {\mathfrak E}_{i} \ar[r]\ar[d] & {\mathfrak E} \ar[r]\ar[d]^{\rm Id} & {\mathfrak E}/{\mathfrak E}_{i} \ar[r]\ar[d] & 0 \\
0 \ar[r] & \tilde{\mathfrak E}_{i}\ar[r]\ar[d] & {\mathfrak E} \ar[r]  & {\mathfrak E}/\tilde{\mathfrak E}_{i}
\ar[r]\ar[d] & 0 \\
         & \tilde{\mathfrak E}_{i}/{\mathfrak E}_{i}\ar[d]  &          &           0                     &   \\
         &                0               &                            &                                 &   \\
}
\end{displaymath}
with ${\mathfrak T}_{i}\cong\tilde{\mathfrak E}_{i}/{\mathfrak E}_{i}$. Since $\tilde{\mathfrak E}_{i}$ and 
${\mathfrak E}_{i}$ are torsion-free, from Section 2 we see that the determinant of 
$\tilde{\mathfrak E}_{i}/{\mathfrak E}_{i}$ is a Higgs bundle, and consequently also is the determinant of 
${\mathfrak T}_{i}$ and we have ${\rm det}\,\tilde{\mathfrak E}_{i}\cong{\rm det}\,{\mathfrak E}_{i}\otimes{\rm det}\,{\mathfrak T}_{i}$. 
If we use this isomorphism and we consider now the saturated flag\footnote{Since ${\mathfrak E}_{i}\subset{\mathfrak E}_{i+1}$, there exists a map 
${\mathfrak E}/{\mathfrak E}_{i}\rightarrow{\mathfrak E}/{\mathfrak E}_{i+1}$ such that the obvious diagram 
commutes. Now, any element in $\tilde{\mathfrak E}_{i}$ can be projected on 
$\tilde{\mathfrak E}_{i}/{\mathfrak E}_{i}\cong{\mathfrak T}_{i}\subset{\mathfrak T}_{i+1}$, so it is zero in 
${\mathfrak E}/\tilde{\mathfrak E}_{i+1}$ and hence $\tilde{\mathfrak E}_{i}\subset\tilde{\mathfrak E}_{i+1}$ and 
$\tilde{\cal F}$ is a flag, which is obviously saturated.} $\tilde{\cal F}=\{(\tilde{\mathfrak E}_{i},n_{
i})\}_{i=1}^{k}$ of ${\mathfrak E}$ we get
\begin{eqnarray*}
 {\mathfrak T}_{\tilde{\cal F}} &=& \prod_{i=0}^{k}(({\rm det}\,\tilde{\mathfrak E}_{i})^{r}\otimes({\rm det}\,{\mathfrak E})^{-r_{i}})^{n_{i}}\\
&\cong& {\mathfrak T}_{{\cal F}'}\otimes\prod_{i=1}^{k}({\rm det}\,{\mathfrak T}_{i})^{rn_{i}}.
\end{eqnarray*}
Since each ${\mathfrak T}_{i}$ is torsion, from a classical result in \cite{Kobayashi} each 
${\rm det}\,{\mathfrak T}_{i}$ admits a nonzero holomorphic section; Now, since $\tilde{\cal F}$ is saturated, 
these conditions on the existence or not of holomorphic sections are satisfied for 
${\mathfrak T}_{\tilde{\cal F}}\otimes{\mathfrak L}$. At this point, by using the isomorphism above if follows 
that the same is true also for the Higgs line bundle ${\mathfrak T}_{{\cal F}'}\otimes{\mathfrak L}$. \;\;Q.E.D.  \\
\begin{pro}
 Let ${\mathfrak E}$ be a torsion-free Higgs sheaf over a compact complex manifold $X$. Then\\
{\bf (i)} If ${\rm rk}\,{\mathfrak E}=1$, then ${\mathfrak E}$ is $T$-stable;\\
{\bf (ii)} If ${\mathfrak L}$ is a Higgs line bundle over $X$, then the tensor product ${\mathfrak E}\otimes{\mathfrak L}$ 
           is $T$-stable (resp. $T$-semistable) if and only if ${\mathfrak E}$ is $T$-stable (resp. $T$-semistable);\\
{\bf (iii)} ${\mathfrak E}$ is $T$-stable (resp. $T$-semistable) if and only if its dual ${\mathfrak E}^{*}$ is 
           $T$-stable (resp. $T$-semistable).
\end{pro}
\noindent {\it Proof:} If ${\rm rk}\,{\mathfrak E}=1$, there are no flags to be considered, hence (i) is trivial. 
Suppose now that ${\mathfrak L}$ is a Higgs line bundle, then in analogy to the classical case, there exists a 
natural correspondence between flags ${\cal F}=\{({\mathfrak E}_{i},n_{i})\}$ of ${\mathfrak E}$ and flags 
${\cal F}\otimes{\mathfrak L}=\{({\mathfrak E}_{i}\otimes{\mathfrak L},n_{i})\}$ of ${\mathfrak E}\otimes{\mathfrak L}$. 
Now, since ${\mathfrak E}$ and ${\mathfrak E}_{i}$ are torsion-free and we have the identities
\begin{equation}
 \bigwedge^{r_{i}}({\mathfrak E}_{i}\otimes{\mathfrak L})\cong(\bigwedge^{r_{i}}{\mathfrak E}_{i})\otimes{\mathfrak L}^{r_{i}}\,, 
 \quad\quad \bigwedge^{r}({\mathfrak E}\otimes{\mathfrak L})\cong(\bigwedge^{r}{\mathfrak E})\otimes{\mathfrak L}^{r} \nonumber 
\end{equation}
we have the following isomorphisms of determinant Higgs bundles 
\begin{equation}
 {\rm det}({\mathfrak E}_{i}\otimes{\mathfrak L})\cong{\rm det}\,{\mathfrak E}_{i}\otimes{\mathfrak L}^{r_{i}}\,, 
 \quad\quad {\rm det}({\mathfrak E}\otimes{\mathfrak L})\cong{\rm det}\,{\mathfrak E}\otimes{\mathfrak L}^r\,. \nonumber
\end{equation}
Now, using these isomorphisms and the expression (\ref{def T}) for the flag ${\cal F}\otimes{\mathfrak L}$ we obtain
\begin{eqnarray*}
 {\mathfrak T}_{{\cal F}\otimes{\mathfrak L}} &=& \prod_{i=1}^{k}\left({\rm det}({\mathfrak E}_{i}\otimes{\mathfrak L})^{r}
                                      \otimes{\rm det}({\mathfrak E}\otimes{\mathfrak L})^{-r_{i}}\right)^{n_{i}} \\
                                  &\cong& \prod_{i=1}^{k}\left(({\rm det}\,{\mathfrak E}_{i})^{r}\otimes
                                          {\mathfrak L}^{r_{i}r}\otimes({\rm det}\,{\mathfrak E})^{-r_{i}}\otimes{\mathfrak L}^{-rr_{i}}\right)^{n_{i}}
                                           \cong {\mathfrak T}_{\cal F} \nonumber
\end{eqnarray*}
and (ii) follows. Finally, assume that ${\mathfrak E}^{*}$ is $T$-stable (resp. $T$-semistable) and let 
${\cal F}=\{({\mathfrak E}_{i},n_{i})\}_{i=1}^{k}$ be a saturated flag of ${\mathfrak E}$. By dualizing the 
Higgs extension of ${\mathfrak E}$ associated to ${\mathfrak E}_{i}$ we get the exact sequence
\begin{displaymath}
 \xymatrix{
0 \ar[r]   &   ({\mathfrak E}/{\mathfrak E}_{i})^{*} \ar[r]  &  {\mathfrak E}^{*} \ar[r]   &  {\mathfrak E}_{i}^{*}
}
\end{displaymath}
with ${\rm rk}({\mathfrak E}/{\mathfrak E}_{i})^{*}=r-r_{i}$ and we obtain from this a flag 
${\cal F}^{*}=\{(({\mathfrak E}/{\mathfrak E}_{i})^{*},n_{i})\}$ of ${\mathfrak E}^{*}$ with 
\begin{equation}
 ({\mathfrak E}/{\mathfrak E}_{k})^{*}\subset\cdots\subset({\mathfrak E}/{\mathfrak E}_{1})^{*}\subset
{\mathfrak E}^{*}\,. \nonumber
\end{equation}
Now, ${\mathfrak E}$ is torsion-free and ${\mathfrak E}/{\mathfrak E}_{i}$ is torsion-free because the flag ${\cal F}$ is saturated, then ${\rm det}\,{\mathfrak E}^{*}\cong({\rm det}\,{\mathfrak E})^{*}$ and ${\rm det}({\mathfrak E}/{\mathfrak E}_{i})^{*}\cong({\rm det}\,{\mathfrak E}/{\mathfrak E}_{i})^{*}$ and hence
\begin{eqnarray*}
 {\mathfrak T}_{{\cal F}^{*}} &=& \prod_{i=1}^{k}\left(({\rm det}({\mathfrak E}/{\mathfrak E}_{i})^{*})^{r}\otimes({\rm det}\,{\mathfrak E}^{*})^{-(r-r_{i})}\right)^{n_{i}} \\
                  &\cong& \prod_{i=1}^{k}\left(({\rm det}\,{\mathfrak E}/{\mathfrak E}_{i})^{-r}\otimes({\rm det}\,{\mathfrak E})^{r-r_{i}}\right)^{n_{i}} \\
                  &\cong& \prod_{i=1}^{k}\left(({\rm det}\,{\mathfrak E}_{i})^{r}\otimes({\rm det}\,{\mathfrak E})^{-r_{i}}\right)^{n_{i}} = {\mathfrak T}_{\cal F}
\end{eqnarray*}
and it follows that ${\mathfrak E}$ is $T$-stable (resp. $T$-semistable). Conversely, assume that 
${\mathfrak E}$ is $T$-stable (resp. $T$-semistable) and let ${\cal F}^{*}=\{({\mathfrak R}_{i},n_{i})\}_{i=1}^{k}$ 
be a saturated flag of ${\mathfrak E}^{*}$. Then for ${\mathfrak R}_{i}$ we have the short exact sequence
\begin{equation}
 \xymatrix{
0 \ar[r]  &  {\mathfrak R}_{i} \ar[r]  &   {\mathfrak E}^{*} \ar[r]  &   {\mathfrak H}_{i} \ar[r] &  0  
} \label{Higgs ext. dual}
\end{equation}
with ${\mathfrak H}_{i}={\mathfrak E}^{*}/{\mathfrak R}_{i}$ torsion-free. By dualizing the Higgs extension 
(\ref{Higgs ext. dual}) we get the exact sequence
\begin{displaymath}
 \xymatrix{
0 \ar[r]   &   {\mathfrak H}_{i}^{*} \ar[r]  &  {\mathfrak E}^{**} \ar[r]   &  {\mathfrak R}_{i}^{*}\,.
}
\end{displaymath}
Since ${\mathfrak E}$ is torsion-free, the natural morphism $\sigma:{\mathfrak E}\rightarrow{\mathfrak E}^{**}$ 
is injective and we can consider the Higgs sheaf ${\mathfrak E}_{i}=\sigma({\mathfrak E})\cap{\mathfrak H}_{i}^{*}$ 
as a Higgs subsheaf of ${\mathfrak E}$ with rank $r_{i}=r-{\rm rk}\,{\mathfrak R}_{i}$. From this we have a flag 
${\cal F}=\{({\mathfrak E}_{i},n_{i})\}$ of ${\mathfrak E}$ with
\begin{equation}
 {\mathfrak E}_{k}\subset{\mathfrak E}_{k-1}\subset\cdots\subset{\mathfrak E}_{1}\subset{\mathfrak E}\,. \nonumber
\end{equation}
Now, we define torsion Higgs sheaves ${\mathfrak T}={\mathfrak E}^{**}/{\mathfrak E}$ and 
${\mathfrak T}_{i}={\mathfrak H}_{i}^{*}/{\mathfrak E}_{i}\subset{\mathfrak T}$. Again, since 
${\mathfrak E}$ is torsion-free, ${\rm det}\,{\mathfrak E}^{**}\cong{\rm det}\,{\mathfrak E}$ and consequently 
${\rm det}\,{\mathfrak T}$ is trivial (as a classical bundle) and 
${\rm det}\,{\mathfrak H}_{i}^{*}\cong{\rm det}\,{\mathfrak E}_{i}$. From this we get
\begin{equation}
 {\rm det}\,{\mathfrak R}_{i} \cong {\rm det}\,{\mathfrak E}^{*}\otimes({\rm det}\,{\mathfrak H}_{i})^{-1}
                             \cong {\rm det}\,{\mathfrak E}^{*}\otimes{\rm det}\,{\mathfrak H}_{i}^{*}
                             \cong {\rm det}\,{\mathfrak E}^{*}\otimes{\rm det}\,{\mathfrak E}_{i}\,. \label{Iso dets}
\end{equation}
Then, from (\ref{Iso dets}) we get
\begin{eqnarray*}
 {\mathfrak T}_{{\cal F}^{*}} &=& \prod_{i=1}^{k}\left(({\rm det}\,{\mathfrak R}_{i})^{r}\otimes({\rm det}\,{\mathfrak E}^{*})^{r_{i}-r}\right)^{n_{i}}\\
                &\cong& \prod_{i=1}^{k}\left(({\rm det}\,{\mathfrak E}_{i})^{r}\otimes({\rm det}\,{\mathfrak E}^{*})^{r_{i}}\right)^{n_{i}}\\
                &\cong& \prod_{i=1}^{k}\left(({\rm det}\,{\mathfrak E}_{i})^{r}\otimes
({\rm det}\,{\mathfrak E})^{-r_{i}}\right)^{n_{i}} = {\mathfrak T}_{\cal F}\\
\end{eqnarray*}
From this isomorphism it follows that ${\mathfrak E}^{*}$ is $T$-stable (resp. $T$-semistable) and hence 
we have proved (iii). \;\;Q.E.D.  \\

\section{The K\"ahler case}
As it is well known \cite{Kobayashi}, if $X$ is K\"ahler there exists a connection between the Mumford-Takemoto 
stability and $T$-stability for coherent sheaves. This result extends naturally to Higgs sheaves and can be written as:
\begin{thm}\label{mu-s-->T-s}
 Let ${\mathfrak E}$ be a torsion-free Higgs sheaf over a compact K\"ahler manifold $X$ with K\"ahler form $\omega$. If ${\mathfrak E}$ is $\omega$-stable (resp. $\omega$-semistable), then it is $T$-stable (resp. $T$-semistable).
\end{thm}
\noindent {\it Proof:} Assume that ${\mathfrak E}$ is $\omega$-stable (resp. $\omega$-semistable) and let ${\cal F}=\{({\mathfrak E}_{i},n_{i})\}_{i=1}^{k}$ be a flag (not necesarily saturated) of ${\mathfrak E}$. Then, as in the classical case we have
\begin{eqnarray*} 
\int_{X}c_{1}({\mathfrak T}_{\cal F})\wedge\omega^{n-1} &=& \sum_{i=1}^{k}n_{i}\int_{X}c_{1}[({\rm det}\,{\mathfrak E}_{i})^{r}\otimes({\rm det}\,{\mathfrak E})^{-r_{i}}]\wedge\omega^{n-1}\\ 
                                            &=& \sum_{i=1}^{k}n_{i}\int_{X}(rc_{1}({\mathfrak E}_{i}) - r_{i}
                                                 c_{1}({\mathfrak E}))\wedge\omega^{n-1} \\
                                            &=& \sum_{i=1}^{k}n_{i}r_{i}r
                                                 (\mu({\mathfrak E}_{i})-\mu({\mathfrak E})) < 0 
\end{eqnarray*}
(resp. $\le 0$). If ${\mathfrak L}$ is a Hermitian flat Higgs line bundle, in particular $c_{1}({\mathfrak L})=0$ 
and we have
\begin{equation}
 {\rm deg}({\mathfrak T}_{\cal F}\otimes{\cal L}) = \int_{X}c_{1}({\mathfrak T}_{\cal F}\otimes{\cal L})\wedge\omega^{n-1} 
                                      = \int_{X}c_{1}({\mathfrak T}_{\cal F})\wedge\omega^{n-1} < 0  \nonumber
\end{equation}
(resp. $\le 0$). Therefore, by using Proposition \ref{Prop. Line bundles} it follows that ${\mathfrak E}$ is 
$T$-stable (resp. $T$-semistable). \;\;Q.E.D.  \\

At this point, as a direct consequence of Theorem \ref{mu-s-->T-s} and Theorem \ref{HK-corresp. sheaves} we 
obtain the following result for reflexive Higgs sheaves over compact K\"ahler manifolds. 
\begin{cor}\label{HYM-->T-ps reflex}
 Let ${\mathfrak E}$ be a reflexive Higgs sheaf over a compact K\"ahler manifold $X$. If ${\mathfrak E}$ has an admissible HYM-metric, then it is $T$-semistable and 
 ${\mathfrak E} = \bigoplus_{i=1}^{s}{\mathfrak E}_{i}$, where each ${\mathfrak E}_{i}$ is a $T$-stable 
 Hermitian-Yang-Mills Higgs sheaf with $\mu({\mathfrak E}_{i})=\mu({\mathfrak E})$. 
\end{cor}

If, on the other hand, we consider locally free Higgs sheaves over compact K\"ahler manifolds, then there exists 
a relation between the notion of $HYM$-metric and the concept of $T$-stability.  In fact, from Theorem \ref{mu-s-->T-s} and 
Theorem \ref{HK-corresp.} we obtain
\begin{cor}\label{HYM-->T-s}
Let ${\mathfrak E}$ be a Higgs bundle over a compact K\"ahler manifold $X$. Then\\
{\bf (i)} If ${\mathfrak E}$ admits a HYM-metric, then it is $T$-semistable and ${\mathfrak E} = \bigoplus_{i=1}^{s}{\mathfrak E}_{i}$, 
where each ${\mathfrak E}_{i}$ is a $T$-stable Hermitian-Yang-Mills Higgs bundle with 
$\mu({\mathfrak E}_{i})=\mu({\mathfrak E})$; \\
{\bf (ii)} If ${\mathfrak E}$ admits an apHYM-metric, then it is $T$-semistable.
\end{cor}

Notice that, the part (i) of Corollary \ref{HYM-->T-s} can be seen as a particular case of Corollary 
\ref{HYM-->T-ps reflex}, however the part (ii) is new. Now, Kobayashi proved in \cite{Kobayashi} a partial 
converse of this Corollary for torsion-free sheaves. The proof of Kobayashi can be easily adapted to Higgs sheaves 
and gives a partial converse of Theorem \ref{mu-s-->T-s}. So we have the following
\begin{thm}\label{T-s-->mu-s, proj case}
 Let ${\mathfrak E}$ be a torsion-free Higgs sheaf over a compact K\"ahler manifold $X$ with K\"ahler form 
$\omega$ and assume either\\
{\bf (a)} The dimension of $H^{1,1}(X,{\mathbb C})$ is equal to one; or\\
{\bf (b)} $\omega$ represents an integral class and ${\rm Pic}(X)/{\rm Pic}^{0}(X)={\mathbb Z}$.\\
If ${\mathfrak E}$ is $T$-stable (resp. $T$-semistable), then it is $\omega$-stable (resp. $\omega$-semistable).
\end{thm}
\noindent {\it Proof:} In analogy to the classical proof of Kobayashi, let ${\mathfrak E}'$ be any Higgs subsheaf 
of ${\mathfrak E}$ with nonzero rank $r'<r$, and consider the Higgs line bundle 
\begin{equation}
 {\mathfrak U}=({\rm det}\,{\mathfrak E}')^{r}\otimes({\rm det}\,{\mathfrak E})^{-r'}\, \label{def U}
\end{equation}
with its degree given by
\begin{equation}
 {\rm deg}\,{\mathfrak U} = \int_{X}(rc_{1}({\mathfrak E}') - r'c_{1}({\mathfrak E}))\wedge\omega^{n-1} 
                          = rr'(\mu({\mathfrak E}') - \mu({\mathfrak E}))\,. \nonumber
\end{equation}
If $[\omega]$ denotes the cohomology class of $\omega$, from either of hypothesis (a) or (b), we obtain 
$c_{1}({\mathfrak U})=a[\omega]$ for some $a\in{\mathbb R}$; hence by integrating this formula it follows that 
${\rm deg}\,{\mathfrak U}=0$ (resp. $>0$) if and only if $a=0$ (resp. $>0$). \\

If ${\mathfrak E}$ is not $\omega$-semistable, there exists a Higgs subsheaf ${\mathfrak E}'$ such that 
$\mu({\mathfrak E}')>\mu({\mathfrak E})$. Then, for the corresponding Higgs line bundle ${\mathfrak U}$ defined by 
(\ref{def U}) we get ${\rm deg}\,{\mathfrak U}>0$ and hence $a>0$. From this we know there exists a positive integer $p$ such that 
${\mathfrak U}^{p}$ admits a nonzero section, say $s$. Now, if $s$ vanishes nowhere on $X$, then 
${\mathfrak U}^{p}$ is trivial as a classical bundle and $c_{1}({\mathfrak U}^{p})=0$ and therefore $a=0$, which 
is a contradiction. Therefore, $s$ necessarily vanishes at some 
point\footnote{Notice that non $T$-semistability means that there exists a nonzero section 
of $T_{\cal F}\otimes{\cal L}$ for some flag ${\cal F}$ and some ${\mathfrak L}$ Hermitian flat Higgs line 
bundle, and such a section vanishes at least at some point.}. In particular, 
${\mathfrak U}^{p}\otimes{\mathfrak L}$ admits a nonzero section, vanishing at some point, for any Hermitian flat 
Higgs line bundle ${\mathfrak L}$ with $L$ trivial. This shows that ${\mathfrak E}$ is not $T$-semistable if it 
is not $\omega$-semistable.\\

If ${\mathfrak E}$ is not $\omega$-stable, then there exists a Higgs subsheaf ${\mathfrak E}'$ such that 
$\mu({\mathfrak E}')\ge\mu({\mathfrak E})$. Now, if the inequality is strict it is also not 
$\omega$-semistable, and hence from the above analysis we conclude that it is not $T$-semistable and in 
particular it is not $T$-stable. If, on the other hand $\mu({\mathfrak E}')=\mu({\mathfrak E})$, it follows that 
${\rm deg}\,{\mathfrak U}=0$ and $a=0$. Hence, the corresponding ${\mathfrak U}$ is flat as a 
classical bundle. Then, by defining $\tilde{\mathfrak L}=(U^{-1},0)$, with $U$ the holomorphic line bundle associated 
to ${\mathfrak U}$, it follows that ${\mathfrak U}\otimes\tilde{\mathfrak L}$ is trivial as a classical bundle 
and hence it admits a nonzero holomorphic section. This shows that ${\mathfrak E}$ is not $T$-stable if it is 
not $\omega$-stable. \;\;Q.E.D.  \\

As it is well known \cite{Hitchin}, there are Higgs bundles over curves that are stable in the sense of 
Mumford-Takemoto, that are not stable as classical bundles. From Theorem \ref{mu-s-->T-s} these bundles 
are $T$-stable as Higgs bundles. Now, from Kobayashi \cite{Kobayashi} it is known that for holomorphic bundles 
over curves, the notions of Mumford-Takemoto stability and $T$-stability are equivalent; hence, such Higgs 
bundles are not $T$-stable in the classical sense. This fact shows that $T$-stability is indeed an extension of 
the classical notion of $T$-stability. Finally, as a direct consequence of Theorem \ref{T-s-->mu-s, proj case} and Theorem \ref{HK-corresp. sheaves} we 
get a partial converse of Corollary \ref{HYM-->T-ps reflex}. To be precise we have the following result
\begin{cor}\label{T-s-->admissible HYM}
 Let ${\mathfrak E}$ be a reflexive Higgs sheaf over a compact K\"ahler manifold $X$ and assume that either (a) 
 or (b) of Theorem \ref{T-s-->mu-s, proj case} holds. If ${\mathfrak E}$ is $T$-stable, then it has an admissible $HYM$-metric.
\end{cor}
This Corollary can be extended to $T$-semistable Higgs sheaves in a very special case. Indeed, if ${\mathfrak E}$ 
is $T$-semistable and ${\mathfrak E} = \bigoplus_{i=1}^{s}{\mathfrak E}_{i}$, where each ${\mathfrak E}_{i}$ is a 
$T$-stable Higgs sheaf with $\mu({\mathfrak E}_{i})=\mu({\mathfrak E})$, then from Corollary 
\ref{T-s-->admissible HYM} each ${\mathfrak E}_{i}$ has an admissible $HYM$-metric and  
(see \cite{Biswas-Schumacher} or \cite{Cardona 2} for details) we get also an admissible $HYM$-metric on 
${\mathfrak E}$. \\

\noindent{\bf Acknowledgements} \\
\noindent This paper was mostly done during a stay of the author at the International School for Advanced Studies 
(SISSA) in Trieste, Italy. The author wants to thank SISSA for the hospitality and support. Finally, the author 
would like to thank U. Bruzzo for some useful comments and suggestions.

\end{document}